\documentclass[10pt]{amsart}

\usepackage{amssymb,eucal}
\usepackage{amscd}
\usepackage[all,cmtip]{xy}
\usepackage{rotating}
\usepackage{hyperref,mathrsfs,soul}
\usepackage{upgreek}
\usepackage{tcolorbox}

\usepackage{amsmath,amssymb}

\setul{}{1.3pt}

\newtheorem{theorem}{Theorem }[section]

\newtheorem{lemma}[theorem]{Lemma}

\newtheorem{remark}[theorem]{Remark}
\newtheorem{corollary}[theorem]{Corollary}

\newcommand{\proj}{\mathrm{proj}}
\def\I{\mathbf{I}}

\def\eop{\hspace*{\fill}{\footnotesize$\blacksquare$}}
\newcommand{\Aut}{\mathrm{Aut}}
\newcommand{\id}{\mathrm{id}}

\newcommand{\mP}{\mathcal{P}}

\newcommand{\mM}{\mathcal{M}}
\newcommand{\mS}{\mathcal{S}}

\newcommand{\mF}{\mathcal{F}}

\newcommand{\bP}{\mathbb{P}}
\newcommand{\mL}{\mathcal{L}}
\newcommand{\mW}{\mathcal{W}}
\newcommand{\F}{\mathbb{F}}
\newcommand{\mH}{\mathcal{H}}

\title[Cyclic STGQs]{Infinite cyclic projective skew translation quadrangles do not exist}

\subjclass[2000]{05B25, 20B10, 20B25, 20E42,
51E12.}

\author{Koen Thas}

\address{{Ghent University},
{Department of Mathematics: Algebra and Geometry},
{Krijgslaan 281, S25, B-9000 Ghent, Belgium}}
\email{koen.thas@gmail.com}
\thanks{}
\date{}

\begin{document}
\maketitle

\begin{abstract}
In this paper we completely classify infinite cyclic projective skew translation quadrangles through a new approach first partially introduced in \cite{Froh}, and in the present paper (unexpectedly) adapted to the infinite case. Very surprisingly, these objects do not exist and only occur in the finite case.   \\
 \end{abstract}

\begin{tcolorbox}
\setcounter{tocdepth}{1}
\tableofcontents
\end{tcolorbox}

\bigskip
\section{Introduction}

There seems to be an unwritten yet strong formalism in incidence geometry which says that (1) usually when some type of object (such as spreads in a projective space of certain dimension) can be constructed in finite geometry, it is much easier to construct the same type of object in the infinite case, and (2) usually classifications and characterizations of objects (such as Moufang quadrangles) tend to be much harder to next-to-impossible in the infinite case. In \cite{BeuCam}, for instance, Beutelspacher and Cameron introduce a 
number of general constructions of objects which live in infinite projective spaces and other infinite geometries, which are hard (and often impossible) to find in the finite case. In \cite{Rosehr}, Rosehr shows that the parameters of infinite lineair translation generalized quadrangles \cite{TGQBook} are not restricted, through a free construction of Kantor familes. In the finite case, there is no such elbow room in construction theory. Also, Tits freely constructed groups with a BN-pair of rank two \cite{POL} which are not related to classical generalized polygons (although they have a huge group of automorphisms!), while in the finite case, the classification of BN-pairs of rank two is known to yield only classical polygons (relying on the classification of finite simple groups) \cite{BUVM}. In the infinite case, classification is out of the question once free constructions come into play. 
On the other hand, in \cite{GSHVM} Steinke and Van Maldeghem manage to obtain a number of infinite analogons to results on finite generalized quadrangles of order $(t,t)$ with opposite axes of symmetry, where the restriction ``$(t,t)$'' is naturally replaced by the fact that the axes of symmetry are assumed to be projective. This approach is the one we will follow in the current paper. 
Steinke and Van Maldeghem are then able to apply their structural results to the case of compact connected quadrangles with topological parameters $(p,p)$ which have opposite axes of symmetry, and show that the quadrangles must be isomorphic to the point-line duals of the symplectic quadrangles $\mW(\mathbb{R})$ or $\mW(\mathbb{C})$. In the finite case, the quadrangles are isomorphic to the point-line duals of the symplectic quadrangles $\mW(t)$. It might be possible that many other (even nonclassical) examples arise in the infinite case once the topological assumptions are dropped altogether. 

Recently, in the paper \cite{Froh}, we have classified cyclic skew translation quadrangles of finite order $(t,t)$, by showing that such quadrangles must be isomorphic to classical symplectic quadrangles $\mW(t)$ with $t$ a prime number. That result turned out to be pivotal for classifying finite skew translation quadrangles of order $(t,t)$ with $t$ even and $t$ not a square in a geometrical fashion \cite{Froh}.
Usually this type of result has a highly specialized proof, and also in this case there seemed not much hope to generalize to the infinite case. On the other hand, although in the finite case only one type of quadrangle emerges, in the aforementioned formalism one would conjecture that (much) more examples arise in the infinite case. It is relevant and definitely insightful to compare the situation with the highly ``acquainted'' problem of classifying cyclic projective planes: in the finite case, it is a major conjecture that such planes are isomorphic to classical planes $\mathbb{P}^2(\F_q)$ defined over a finite field $\F_q$, while in the infinite case, it is easy to freely construct (nonclassical) planes with Singer groups isomorphic to the infinite cyclic group $(\mathbb{Z},+)$ \cite{Hall,Karzel,HP}; by Karzel \cite{Karzel} no infinite classical plane $\mathbb{P}^2(k)$ can admit $(\mathbb{Z},+)$ as a Singer group (but many infinite classical planes {\em do} admit noncyclic Singer groups on the other hand \cite{KTJNT}). Our main result will be even more dramatic.  

In the present paper, we completely classify infinite cyclic projective skew translation quadrangles, and this result comes with two major surprises: first of all, we manage to use an idea of \cite{Froh}, in combination with the theory of divisible groups and a result only recently obtained in 2019 \cite{TGQinf}, to produce a full classification; secondly, 
the outcome of the theorem is that {\em no such generalized quadangles exist}, contrary to what our formalism says, and in stark contrast with for instance Rosher's result \cite{Rosehr} in which he freely constructs certain Kantor families in {\em abelian} groups! As there are next-to-no definite classification results known for infinite generalized quadrangles  with only local group actions, it makes the result the more remarkable. Using a new approach to construct root-elations (in fact, {\em symmetries}) through our theory of angles introduced in \cite{STGQ,Froh} plus further ingredients from the theory of divisible groups and infinite translation quadrangles \cite{TGQinf}, we are able to control the elements of the associated Kantor family in a precise manner, which is a big surprise to the author. Maybe even more striking than the main result itself, is the fact that it is possible to obtain it in such a concise, geometrical way (given the pain it takes in the finite case \cite{Froh}). 

We feel that our new approach to construct and control root-elations will serve as a springboard to similar constructions in a wide variety of problems in the infinite theory of generalized quadrangles which require existence of root-elations. Especially in the local case such results are extremely sparse.  

As a by-product, we obtain a complete classification of projective cyclic skew translation quadrangles (Corollary \ref{corcyc}).\\

\subsection*{Structure of the paper}

In the next two sections (\S \ref{intro} and \S \ref{KF}) we introduce some notions needed to comprehend this paper in a self-contained way. In sections \ref{proj} and \ref{notes}, we make a number of initial observations which prepare the proof of the main theorem. In the final section \ref{main}, the main result is obtained.



\bigskip
\section{Combinatorics}
\label{intro}

Let $\Upgamma$ be a thick {\em generalized quadrangle} (GQ). It is a rank $2$ geometry $\Upgamma = (\mP,\mL,\I)$ (where we call the elements of $\mP$ ``points'' and those of $\mL$ ``lines'')
such that the following axioms are satisfied:
\begin{itemize}
\item[(a)]
there are no ordinary digons and  triangles contained in $\Upgamma$;
\item[(b)]
each two elements of $\mP \cup \mL$ are contained in an ordinary quadrangle;
\item[(c)]
there exists an ordinary pentagon.
\end{itemize}
It can be shown that there are exist constants $s$ and $t$ such that each point is incident with $t + 1$ lines and each line is incident with $s + 1$ points. We say that
$(s,t)$ is the {\em order} of $\Upgamma$. 
Note that an ordinary quadrangle is just a (necessarily thin) GQ of order $(1,1)$ | we call such a subgeometry also ``apartment'' (of $\Upgamma$).

\subsection{Subquadrangles}

A {\em subquadrangle} (subGQ) $\Upgamma' = (\mP',\mL',\I')$ of a generalized quadrangle $\Upgamma = (\mP,\mL,\I)$, is a GQ for which 
$\mP' \subseteq \mP$, $\mL' \subseteq \mL$, and $\I'$ is the incidence relation which is induced by $\I$ on 
$(\mP' \times \mL') \cup (\mL' \times \mP')$. We say that $\Upgamma'$ is {\em ideal} if for each of its points $x$, any line of $\Upgamma$ incident with $x$ is also 
a line of $\Upgamma'$.

\subsection{Regularity}
\label{reg}

Let $\Upgamma$ be a thick GQ of order $(s,t)$. If $X$ and $Y$ are (not necessarily different) points, or lines, $X \sim Y$ denotes 
the fact that there is a line incident with both, or a point incident with both. Then $X^{\perp} := \{ Y \ \vert\ Y \sim X \}$, and 
if $S$ is a point set, or a line set, $S^{\perp} := \cap_{s \in S}s^{\perp}$ and $S^{\perp\perp} := {(S^{\perp})}^{\perp}$. 
A particularly important example is the case where $S = \{X,Y\}$ is a set of distinct noncollinear points (or nonconcurrent lines, but this is 
merely the dual situation, which we leave to the reader); then each line incident with $X$ is incident with precisely one point of $\{ X,Y\}^{\perp}$ (so if $\Upgamma$ is finite, 
$\vert \{ X,Y\}^{\perp} \vert = t + 1$). The set $\{ X,Y\}^{\perp\perp}$ consists of all points which are collinear with every point of $\{X,Y\}^{\perp}$, 
so 
\begin{equation}
\{ X, Y\} \subseteq \{X,Y\}^{\perp\perp}. 
\end{equation}

Let $Z$ be any point of $\{X,Y\}^{\perp}$; if each line incident with $Z$ is incident with exactly one point of $\{X,Y\}^{\perp\perp}$, then 
this property is independent of the choice of $Z$, and we say that $\{ X,Y\}$ is {\em regular}. (In the finite case, we could equivalently have asked 
that $\vert \{ X,Y\}^{\perp\perp} \vert = t + 1$.)  We call a point/line $X$ {\em regular} if $\{ X,Y \}$ is regular for all $Y \not\sim X$.

\section{EGQs, STGQs and Kantor families}
\label{KF}

\subsection{Symmetry}

Isomorphisms and automorphisms of generalized quadrangles are defined in the usual manner. See chapter 1 of \cite{PT}. The automorphism group of a 
GQ $\Upgamma$ will be denoted by $\Aut(\Upgamma)$. 

Let $X$ be a point or a line in a thick GQ $\Upgamma$. A {\em symmetry} with {\em center} $X$ (in the case of a point) or {\em axis} $X$ (in the case of a line)
is an element of $\Aut(\Upgamma)$ that fixes each element of $X^{\perp}$. We say that $X$ is a {\em center of symmetry} (point case) or an {\em axis of symmetry} (line 
case) if there exist  $Y$ and $Z$  in $X^{\perp}$ such that $Y \not\sim Z$, for which the group of all symmetries $\mS(X)$ with center/axis $X$ acts 
transitively on $\{ Y,Z\}^{\perp} \setminus \{X\}$. This definition does not depend on the choice of $(Y,Z)$, and one easily shows that ``transitive'' implies
``sharply transitive.'' Note that a center/axis of symmetry is necessarily regular.




\subsection{Elation generalized quadrangles}

Let $\Upgamma$ be a generalized quadrangle and let $x$ be one of its points. If there is an automorphism group $E$ which fixes each line incident with $x$ and acts 
sharply transitively on the set of points which are not collinear with $x$, then we say that $\Upgamma$ is an {\em elation generalized quadrangle} (EGQ) with {\em elation group} $E$  and {\em elation point} $x$. We often write $(\Upgamma, E)$ if we want to make the elation group explicit. For more on elation generalized quadrangles, we refer to the monograph \cite{LEGQ}.

\subsection{Skew translation generalized quadrangles}

If furthermore $x$ is a center of symmetry for which the associated group of symmetries $\mS(x)$ is a subgroup of $E$, then we say that 
$\Upgamma$ is a {\em skew translation generalized quadrangle} (STGQ). STGQs are very common geometries. Each known finite generalized quadrangle which is not isomorphic to the Hermitian quadrangle $\mH(4,q^2)$ or its point-line dual for some finite field $\F_q$, 
is | modulo point-line duality | either an STGQ, or the Payne-derived quadrangle of an STGQ. (More details on this statement can be found in \cite{STGQ}.) 

In the infinite case not much is known about STGQs, and one of the goals of this paper is to contribute to this caveat. 

\subsection{Kantor families}

Let $(\Upgamma, E)$ be an EGQ with elation point $x$. Let $y$ be any point not collinear with $x$. For any line $U \I y$, let $[U]$ be the unique line incident with $x$ which meets $U$; the map 
\begin{equation}
\Big[\ \cdot\ \Big]\ : \ V \ \mapsto\ [V] 
\end{equation}
defines a bijection between the lines incident with $y$ and the lines incident with $x$. For any $V \I y$, define $E_V^*$ as $E_{V \cap [V]}$, the point-stabilizer 
of $V \cap [V]$ in $E$. Then $(\mF, \mF^*)$, with $\mF = \{ K_U\ \vert\ U \I y \}$ and $\mF^* = \{ K_U^*\ \vert\ U \I y \}$, is called the {\em Kantor family} associated to the EGQ $(\Upgamma, E)$. (The choice of the point $y$ is not important here, as $E$ acts transitively on the points not collinear with $x$.)  Vice versa, Kantor families can be defined axiomatically (in a purely group-theoretical fashion), and give rise to EGQs \cite{LEGQ}.  

In terms of Kantor families, a Kantor family 
$(\mF,\mF^*)$ gives rise to an STGQ if and only if there is a normal subgroup $\mS$ of $E$ such that for each $A \in \mF$, we have $A^* = A\mS$ (see \cite{STGQbirth} or \cite[chapter 10]{PT});  $\mS$ then corresponds to the group of symmetries with center the elation point. Note that for different $A, B$ in $\mF$, we have $A^* \cap B^* = \mS$.



\section{The projective plane $\Upgamma_x$}
\label{proj}

Let $x$ be a regular point of the generalized quadrangle $\Upgamma$. Define a new rank $2$ geometry $\Upgamma_x$  as follows:
\begin{itemize}
\item[(A)]
the \ul{points} are the lines incident with $x$, and the sets $\{ u, x \}^{\perp}$ with $u$ not collinear with $x$;
\item[(B)]
the \ul{lines} are the points collinear with $x$ (including $x$);
\item[(C)]
the \ul{incidence} is (inverse) containment. 
\end{itemize} 

If $\Upgamma$ is finite and its order is $(t, t)$, then one easily shows that $\Upgamma_x$ is a projective plane of order $t$. In the infinite case, asking that the order 
is $(t, t)$ is not sufficient: we say that $x$ is {\em projective}, or that $\Upgamma$ is {\em projective} at $x$ if $\Upgamma_x$ is a projective plane. Asking that $x$ is projective is a natural (and very good) way to express that the number of lines incident with a point equals the number of points incident with a line.

\section{Note on structure of the elements in $\mF$}
\label{notes}

In this section we suppose that $(\Upgamma, K)$ is an infinite projective STGQ with elation point $x$, and so that the group of symmetries $\mS = \mS(x)$ is cyclic (hence isomorphic to $(\mathbb{Z},+)$). The associated Kantor family is $(\mF,\mF')$. 

Since $x$ is projective, we know that $\Upgamma_x$ is a projective plane. Now $K$ naturally induces an automorphism group of $\Upgamma_x$, and obviously 
the kernel of this action is $\mS$. Since $K/\mS$ fixes the line $x$ of the plane $\Upgamma_x$ and also all its points (which corresponds to the quadrangle lines incident with $x$),  one readily sees that $\Upgamma_x$ is a projective translation plane with translation group $T := K/\mS$ and translation line $x$.  
So we know (see for instance \cite{Knarr}) that $K/\mS$ is isomorphic to the additive group of a left or right vector space $V$ over some division ring $\wp$. Let $(T,\{ T_i \}_i)$ be the corresponding congruence partition \cite[chapter VII, \S 3]{HP}, and let $A \in \mF$ be arbitrary. Then $A \cong A\mS/\mS$ is isomorphic to a member of $\{ T_i \}_i$, and so $A$ is isomorphic to the additive group of a left or right vector space over $\wp$ as well. Two ingredients arise which we will need below:
\begin{itemize}
\item
If $\wp$ has characteristic $0$, any additive group $(G, +)$ of a vector space over $\wp$ is a {\em divisible group}: the equation $X^N =  g$, where we use exponential notation and with $g \in G$ and 
$N$ any positive integer different from $0$, always has a solution. (This is not true in positive characteristic.) It is easy to see this: the center of $\wp$ is a field of characteristic $0$, so it contains $\mathbb{Q}$. So $G$ is a $\mathbb{Q}$-vector space.  
\item
If $\wp$ has characteristic $p > 0$, $p$ a prime, then each element $a$ of $G$ satisfies the equality $a^p = \id$. 
\end{itemize}

\section{Proof of the main result}
\label{main}

In this section we suppose that $(\Upgamma, K)$ is an infinite projective STGQ with elation point $x$, and so that the group of symmetries $\mS = \mS(x)$ is cyclic (hence isomorphic to $(\mathbb{Z},+)$). The associated Kantor family is $(\mF,\mF')$. Our goal is to show that such $\Upgamma$ cannot exist.

\subsection{General considerations}

Let $A$ be any member of $\mF$, and let it act on $\mS$ by conjugation. As the automorphism group of $\mS \cong \Big(\mathbb{Z},+\Big)$ is isomorphic to 
$C_2$ (the only nontrivial automorphism being an involution determined by $\upgamma: 1 \mapsto -1$), it follows that the action of $A$ on $\mS$ has a kernel 
$N(A)$ of at most index $2$ in $A$. If for each $A \in \mF$ we would have that $N(A) = A$, then $[A,\mS] = \{ \id \}$ and since the elements of $\mF$ generate 
the elation group $K$, it follows that $\mS \leq Z(K)$.

\begin{remark}{\rm 
Let $U, V$ be nonconcurrent lines which both meet $[V] \I x$, and suppose that $\eta \in N(K_U)$, while $\eta \in K_V \setminus N(K_V)$. As $\eta \in N(K_U)$, we have that $\eta$ commutes with $\mS$. On the other hand, as $\eta$ fixes $V$ and commutes  with $\mS$, it follows that $\eta \in N(K_V)$, contradiction. So if $\eta \in N(K_U)$, then either it fixes \ul{all lines on $[U] \cap V$, or only $[U]$ on $[U] \cap V$}.}\\
\end{remark}   

\subsection{Angles}

Let $(\Upgamma, K)$ be an STGQ with elation point $x$, and associated group of symmetries (with center $x$) $\mS$. Let $(\mF,\mF^*)$ be the associated Kantor family, and let $A \in \mF$ be arbitrary. Let $A$ fix the line $U \sim [U] \I x$, where $U$ is not incident with $x$. Suppose that $A$ fixes each point incident with $[U]$. Let $V \sim [U]$, with $V$ not incident with $x$; then the {\em angle} of $a$ at $V$ is the unique element $s$ of $\mS$ for which $V^a = V^s$.

\subsection{Characteristic 0}

Let $a \in N(A)^\times$, and note that $a$ commutes with $\mS$ by definition of $N(A)$. 

We know that $A$ is isomorphic to the additive group of a left or right vector space over some division ring $\wp$, and in this section we suppose that $\wp$ has characteristic $0$, so 
$A$ is divisible. A quotient of a divisible group is also divisible, so if $N(A)$ would have index $2$, $A/N(A)$ would be a finite group which is divisible | which is only possible if it is trivial. So $N(A) = A$ and hence $A$ commutes with $\mS$. As we have seen, it follows that $\mS$ is a subgroup of the center of $K$. Let $A \in \mF$ be arbitrary, and let $a \in A^\times$. Suppose $A = K_U$ with $U$ a line not incident with $x$, and consider an arbitrary point $y \ne x$ incident with $[U]$. As $\mS$ commutes with $A$, it follows that either $[U]$ is the only fixed line of $a$ incident with $y$, or all lines incident with $y$ are fixed. This readily implies that per point of $[U]$, at most one nontrivial angle occurs for any element of $A$. In other words: let $b \in A^\times$ and $w \I [U]$. Then either all lines on $w$ are fixed by $b$, or $b$ acts as some element of $\mS$ on the lines incident with $w$. 
Take an element $a \in A^\times$ which is not a symmetry with axis $[U]$, and consider a nontrivial angle $\upgamma$ on some point $u \I [U]$. Choose a positive integer $N$ such that no element $\alpha$ exists in $\mS$ for which $\alpha^N = \upgamma$. In $A$, an element $b$ exists such that $b^N = a$ as $A$ is divisible. This means that on $u$, we have an angle $\alpha$ for which $\alpha^N = \upgamma$, contradiction. It follows that all elements of $A$ are symmetries with axis $[U]$, and since $A$ was arbitrary, this means that \ul{$K$ is abelian and $\Upgamma$ is a TGQ}. \\

By \cite{TGQinf}, it follows that $(\Upgamma, K)$ is an ideal subGQ of a linear TGQ $(\widetilde{\Upgamma}, \widetilde{K})$ (which is a TGQ which allows a projective representation over a division ring $\Game$), for which $K \leq \widetilde{K}$ (and which has the same translation point $x$). Since $\Upgamma$ is ideal in $\widetilde{\Upgamma}$ and since $\mS \leq \widetilde{K}$, we have that the translation point $x$ is also a center of symmetry of $\widetilde{\Upgamma}$ with associated group $\mS$. 
Since $\widetilde{\Upgamma}$ is linear, it follows that $\mS$ is isomorphic to the additive group of a left or right vector space over $\Game$. But that contradicts the fact that $\mS$ is cyclic. \eop \\


 

\subsection{Characteristic $p \ne 0$}


Let $a \in N(A)^\times$ and let $\upgamma$ be a nontrivial angle of $a$. Then $a^{-1}\upgamma$  has order $p$, and hence $\upgamma^p = \id$ (since $[a,\upgamma] = \id$). As there are no nontrivial elements with finite order in $\Big(\mathbb{Z},+\Big)$, it follows that there are no nontrivial angles, which means that $a$ is a line symmetry.  \\

Suppose that for some $B \in \mF$ we have that $N(B) \ne B$: then $B/N(B) \cong C_2$. It follows that $p = 2$. For, suppose $p$ is odd and that $B/N(B) \cong C_2$. Then for all $\beta \in B$ we must have that $\beta^2 \in N(B)$. But as $p$ is odd, we have that $\langle \beta^2 \rangle = \langle \beta \rangle \leq N(B)$ for all $\beta \in B$, so that $N(B) = B$, contradiction.\\

So if $p > 2$, we already know that $N(A) = A$ for all $A \in \mF$, and that all elements in $A$ are line symmetries. So again, \ul{$K$ is abelian and $\Upgamma$ is a TGQ}. Now suppose $p = 2$, and let $B \ne N(B)$ for some $B \in \mF$. Then $B/N(B) \cong C_2$. Let $B$ fix the line $V$ not incident with $x$. Let $v \ne x$ be any point incident with $[V]$; then for $b \in B$ we have that either it fixes all lines on $v$, or it fixes at most one line incident with $v$ and different from $[V]$. To see this, identify the action of $b$ on the lines incident with $v$ and different from $[U]$ with the action of $b$ on $\mS$, and suppose that $b \not\in N(B)$. To find a contradiction, suppose that $b$ fixes two different lines $V, W$ incident with  $v$ and different from $[V]$. Let $\upgamma \in \mS$ be the point-symmetry which sends $V$ to $W$; then $b$ and $\upgamma$ commute, and so $[\upgamma,b] = \id$. From 
\begin{equation}
b\upgamma b \ =\ \upgamma
\end{equation}   
(while noting that $b$ is an involution) follows that the nontrivial automorphism $b$ of $\mS \cong (\mathbb{Z},+)$ fixes $\upgamma$. But the only nontrivial automorphism (generated by $\alpha: 1 \mapsto -1$) of the latter only fixes $\id$. It follows that $b$ only fixes $[U]$ and $V$. \\


Next, we show that such elements cannot exist. Let $A \in \mF$ and suppose $N(A) \ne A$, so that $N(A)$ has index $2$ in $A$. Suppose that $A$ fixes the line $U$ which is not incident with $x$, and let $a \in A \setminus N(A)$. Note that for each line $U'$ which meets $[U]$ and which is not incident with $x$, we also have that $N(K_{U'}) \ne K_{U'}$ (as $K$ acts transitively on these lines). Let $B$ and $C$ be different elements in $\mF$; then $N(B)$ has at most index $2$ in $B$ and the same holds for $[C : N(C)]$. It is easy to see that 
there exists an $N(B)$-orbit $\Omega(B)$ in $[U] \setminus \{x\}$ and an $N(C)$-orbit $\Omega(C)$ in $[U] \setminus \{x\}$ which have at least $2$ points in common. Suppose $u$ and $v$ are two such (different) points. without loss of generality, we suppose that $U \I u$. Now let $\upgamma \in N(B)$ send $u$ to $v$, and let $\upgamma' \in N(C)$ also send $u$ to $v$. As $\upgamma$ and $\upgamma'$ are both symmetries, they both commute with $a$; for instance, the commutator $[a,\upgamma]$ is both a line-symmetry in $B$ and 
a point-symmetry in $\mS$, while $\mS \cap B = \{ \id \}$. It follows that $a$ fixes both lines $U^\upgamma$ and $U^{\upgamma'}$, which are different lines incident with $v$. This is a contradiction. 

It follows that for all primes $p$ and all elements $A \in \mF$, we have that $A = N(A)$, so that as above, $\Upgamma$ is a TGQ and $K$ is abelian. By \cite{TGQinf}, we know that $K$ is isomorphic to the additive group of an $\mathbb{F}_p$-vector space, implying in particular that each element of $\mS$ has order $p$. As $\mS \cong (\mathbb{Z},+)$, this is our final contradiction. \eop \\

Together with the classification obtained in \cite{Froh}, we obtain our final result:

\begin{corollary}[Classification of projective cyclic STGQs]
\label{corcyc}
Let $(\Upgamma, K)$ be a projective cyclic STGQ. Then $\Upgamma$ is finite and isomorphic to a symplectic quadrangle $\mW(p)$ for some prime $p$, while 
$K$ is an elementary abelian $p$-group. 
\eop \\
\end{corollary}

\begin{corollary}[Implication for Kantor families]
Let $(\mF,\mF^*)$ be a Kantor family in the group $K$, such that for each group $A \in \mF$ we can write $A^* = A\mS$ for the normal subgroup $\mS$ of $K$. If $\mS$ is cyclic, then 
$K$ is finite and isomorphic to an elementary abelian $p$-group $C_p \times C_p \times C_p$ for some prime $p$.  
\eop \\
\end{corollary}

\end{document}